\documentclass[12pt,a4paper,reqno]{amsart}
\usepackage{amsmath}
\usepackage{amsfonts}
\usepackage{amssymb}
\numberwithin{equation}{section}

\theoremstyle{plain}

\newtheorem{theorem}[subsection]{Theorem}
\newtheorem{proposition}[subsection]{Proposition}
\newtheorem{lemma}[subsection]{Lemma}

\theoremstyle{definition}

\newtheorem{question}[subsection]{Question}

\renewcommand{\leq}{\leqslant}
\renewcommand{\geq}{\geqslant}

\newsavebox{\proofbox}
\savebox{\proofbox}{\begin{picture}(7,7)%
  \put(0,0){\framebox(7,7){}}\end{picture}}

\newcommand{\md}[1]{\ensuremath{(\mbox{mod}\, #1)}}

\def\proof{\noindent\textit{Proof. }}
\def\endproof{\hfill{\usebox{\proofbox}}}
\def\E{\mathbb{E}}
\def\Z{\mathbb{Z}}
\def\R{\mathbb{R}}
\def\T{\mathbb{T}}
\def\C{\mathbb{C}}

\def\re{\mbox{Re}}

\def\ni{\noindent}
\def\vs{\vspace{11pt}}

\begin{document}

\title{Freiman's theorem in an arbitrary abelian group}

\author{Ben Green}
\address{Department of Mathematics, University of Bristol, University Walk, Bristol BS8 1TW, England
}
\email{b.j.green@bristol.ac.uk}

\author{Imre Z.~Ruzsa}
\address{Alfr\'ed R\'enyi Mathematical Institute, Hungarian Academy of Sciences, Budapest, Pf. 127, H-1364 Hungary
}
\email{ruzsa@renyi.hu}

\thanks{While this work was carried out, the first author was supported by a PIMS postdoctoral fellowship at the University of British Columbia, Vancouver, Canada.}

\begin{abstract} 
A famous result of Freiman describes the structure of finite sets $A \subseteq \mathbb{Z}$ with small doubling property. If $|A + A| \leq K|A|$ then $A$ is contained within a multidimensional arithmetic progression of dimension $d(K)$ and size $f(K)|A|$. Here we prove an analogous statement valid for subsets of an arbitrary abelian group.
\end{abstract}

\maketitle

\section{Introduction}\label{sec1} \ni
Throughout this paper $A$ will be a finite subset of a (not necessarily finite) abelian group $G$. We define the sumset $A + A$ to be the collection of all sums $a + a'$, $a,a' \in A$. If $|A + A| = K|A|$ then we say that $A$ has \textit{doubling} $K$, and if $K$ is very small in comparison to $|A|$ (we will be deliberately vague about what we mean by this) then we say that $A$ has \textit{small doubling}.\vs

\ni Freiman's theorem \cite{freiman} gives a description of sets of integers with small doubling. If $A \subseteq \mathbb{Z}$, and if $|A + A| \leq K|A|$, then $A$ is contained in a proper arithmetic progression of dimension $d(K)$ and size at most $f(K)|A|$. Recall that an arithmetic progression of dimension $d$ and size $L$ is a set of the form
\begin{equation}\label{eq1} P \; = \;  \left\{ v_0 + l_1 v_1 + \dots + l_d x_d \; : \; 0 \leq l_j < L_j\right\},\end{equation}
where $l_1l_2\dots l_d = L$. $P$ is said to be proper if all of the sums in \eqref{eq1} are distinct, in which case $|P| = L$. Observe that Freiman's theorem gives, ignoring quantitative issues connected with the dependence on $K$, a complete description of sets with small doubling. Indeed if $A$ is contained in a proper progression of dimension $d(K)$ and size $f(K)|A|$ then it is easy to see that $|A + A| \leq 2^{d(K)}f(K)|A|$. \vs

\ni Our aim in this paper is to give a similarly complete description of sets with small doubling in an arbitrary abelian group $G$.\vs

\ni The notion of a proper progression makes perfect sense in this more general setting, but the most na\"{\i}ve attempt to generalise Freiman's theorem meets with an obvious failure. If $A = \mathbb{F}_2^k$ then $|A + A| = |A|$, and so $A$ has extremely small doubling. However $A$ is not contained in a progression with small dimension.\vs

\ni In general, any large subset of the product of a group and a proper progression will have small doubling. By a \textit{coset progression} of dimension $d$ we will mean a subset of $G$ of the form $P + H$, where $H \leq G$ is a subgroup, $P$ is a proper progression of dimension $d$ and the sum is direct in the sense that $p + h = p' + h'$ if and only if $h = h'$ and $p = p'$. By the \textit{size} of a coset progression we mean simply its cardinality.

\begin{theorem}\label{mainthm} Let $A \subseteq G$ satisfy $|A + A| \leq K|A|$. Then $A$ is contained in a coset progression of dimension at most $d(K)$ and size at most $f(K)|A|$. We may take $d(K) = CK^4 \log (K+2)$ and $f(K) = \exp(CK^4 \log^2 (K+2))$ for some absolute constant $C$.
\end{theorem}

\ni This result clearly implies Freiman's theorem for subsets of $\mathbb{Z}$.\vs 

\ni The only other types of group for which a complete Freiman type result was previously known are the groups with torsion bounded by $r$. In this case a theorem of the second author \cite{ruzsa-freiman-torsion} asserts that if $G$ has torsion bounded by $r$ and if $A \subseteq G$ satisfies $|A + A| \leq K|A|$, then $A$ is contained in a coset of a subgroup of cardinality at most $K^2 r^{K^4}|A|$. These bounds were improved somewhat in \cite{green-ruzsa}; see also \cite{dhp} for more detailed information in the case $G = (\Z/2\Z)^n$, $K < 4$.\vs
 
\ni Theorem \ref{mainthm} also implies a result of this type, albeit with worse bounds. Indeed if $G$ has torsion bounded by $r$ and if $S \subseteq G$ is a grid of dimension $d$, then $S$ is actually contained a coset of a subgroup of size at most $r^d|S|$.\vs

\ni A different proof of Freiman's theorem was given by the second author \cite{ruzsa-freiman}. To obtain Theorem \ref{mainthm} we will follow the broad scheme of this proof. We will also incorporate some important refinements due to Chang \cite{chang}.\vs

\ni A complete account of the arguments of \cite{chang} and \cite{ruzsa-freiman}, together with all the necessary background, may be found in lecture notes of the first author \cite{green-edinburgh}. We will also draw on results from \cite{green-ruzsa}, and for these reasons this paper is far from self-contained. To conclude this introduction we will however recall the broad outline of \cite{ruzsa-freiman}, indicating the refinements of Chang, and remarking on what further modifications are required to prove Theorem \ref{mainthm}. \vs

\ni Before doing this let us recall two important pieces of nomenclature. Let $s \geq 2$ be an integer, let $G,G'$ be two abelian groups and let $A \subseteq G$ and $A' \subseteq G'$ be sets. Let $\phi : A \rightarrow A'$ be a map. Then we say that $\phi$ is a \textit{Freiman $s$-homomorphism} if whenever $a_1,\dots,a_s,b_1,\dots,b_s \in A$ satisfy
\[ a_1 + a_2 + \dots + a_s \; = \; b_1 + b_2 + \dots + b_s\]
we have
\[ \phi(a_1) + \phi(a_2) + \dots + \phi(a_s) \; = \; \phi(b_1) + \phi(b_2) + \dots + \phi(b_s).\]
If $\phi$ has an inverse which is also an $s$-homomorphism then we say that $\phi$ is a Freiman $s$-isomorphism, and write $A \cong_s A'$. Later on, we will use the the fact that the image of a proper coset progression under a $2$-isomorphism is another proper coset progression with the same dimension and size; we leave the proof as an exercise to the reader.\vs

\ni Now let $\Gamma \subseteq \widehat{G}$ be a set of $d$ characters $\gamma_1,\dots,\gamma_d$ on an abelian group $G$. If $\rho > 0$ is a real number then we define the \textit{Bohr set} $B(\Gamma,\rho)$ to consist of all $x \in G$ for which $(2\pi)^{-1}|\arg(\gamma_j(x))| \leq \rho$ for all $j = 1,\dots,d$. We refer to the parameter $d$ as the \textit{dimension} of the Bohr set.\vs

\ni Now for the (very rough) outline of \cite{ruzsa-freiman}. Suppose that $A \subseteq \mathbb{Z}$ has doubling constant $K$.\vs

\ni\textbf{Step 1} (Finding a good model). There is a prime $p$, $K^{16}|A| \leq p \leq  2K^{16}|A|$, and a set $A' \subseteq A$, $|A'| \geq |A|/2$, such that $A'$ is 8-isomorphic to a subset of $\mathbb{Z}/p\mathbb{Z}$.\vs

\ni\textbf{Step 2} (Bogolyubov's argument). If $X$ is a large subset of $\mathbb{Z}/p\mathbb{Z}$ then $2X - 2X$ contains a large Bohr set.\vs

\ni\textbf{Step 3} (Structure of Bohr sets). A large Bohr set in $\mathbb{Z}/p\mathbb{Z}$ contains a large proper progression with small dimension.\vs

\ni\textbf{Step 4} (Pullback). The structure obtained in step 3 may be used to find a progression of small dimension in which the original set $A$ is contained.\vs

\ni The refinements of Chang allow one to find a larger Bohr set in step 2, and to perform the pullback operation of step 4 more economically.\vs

\ni In order to obtain Theorem \ref{mainthm} we must generalise each of these four steps to general abelian groups. Modifying steps 2 and 4 is straightforward. In step 3, we must elucidate the structure of a Bohr set in a general abelian group. This is not particularly difficult and is accomplished in \S \ref{sec4}. The most interesting part of our argument is the appropriate generalisation of step 1. The form of this generalisation, which clarifies our use of the term \textit{model} above, is the following:
\begin{proposition}[Good models]\label{prop2} Let $A \subseteq G$ be a set with doubling constant $K$. Let $s \geq 2$ be an integer. Then there is a group $G'$, $|G'| \leq (10sK)^{10K^2}|A|$, such that $A$ is Freiman $s$-isomorphic to a subset of $G'$.
\end{proposition}
\ni The proof of this result is the objective of \S \ref{sec2}. Its purpose is that it allows one to work, up to Freiman isomorphism, in an environment where the tools of harmonic analysis may be used efficiently. It is not completely obvious that a given set $A$ has a \textit{finite} model. We will need this fact to prove Proposition \ref{prop2}, and we will prove it at the start of \S \ref{sec2}.\vs

\ni\textit{Notation and Tools.} At various points in the paper it will be convenient to work with the \textit{normalized} counting measure on $G$, since this avoids the appearance of $|G|$ in our formul{\ae}. We use the notation of expectation: if $f : G \rightarrow \C$ is a function then $\E f = \E_{x \in G} f(x) := |G|^{-1}\sum_{x \in G} f(x)$. If $\gamma \in \widehat{G}$ is a character then we write
\[ \widehat{f}(\gamma) = f^{\wedge}(\gamma) := \E_{x \in G} f(x)\gamma(x)\]
for the Fourier transform of $f$ at $\gamma$. We have \textit{Plancherel's identity}
\[ \E_{x \in G} f(x) \overline{g(x)} = \sum_{\gamma \in \widehat{G}} \widehat{f}(\gamma) \overline{\widehat{g}(\gamma)}\]
and the \textit{inversion formula}
\[ f(x) = \sum_{\gamma \in \widehat{G}} \widehat{f}(\gamma) \overline{\gamma(x)}.\]
Defining the \textit{convolution} of two function $f$ and $g$ by
\[ (f \ast g)(x) := \E_{y \in G} f(y) g(x - y),\]
we see that taking Fourier transforms converts convolution into multiplication, that is to say $(f \ast g)^{\wedge} = \widehat{f}\widehat{g}$. Finally, if $S$ is a set then we write $1_S$ for the characteristic function of $S$, defined by $1_S(x) = 1$ if $x \in S$ and $1_S(x) = 0$ otherwise. Note that if $S \subseteq G$ then $\E 1_S = |S|/|G|$; this quantity should be thought of as the relative density of $S$ in $G$.\vs

\ni At various points in the sequel we will use an inequality of the second author \cite{ruzsa-plun}, generalizing and simplifying earlier inequalities of Pl\"unnecke \cite{plunnecke}. If $G$ is an abelian group, and if $A \subseteq G$ is a set, then for integers $k,l \geq 0$ we define $kA - lA$ to be the collection of all sums $a_1 + \dots + a_k - a'_1 - \dots - a'_l$, where $a_i, a'_i \in A$.

\begin{proposition}[\cite{ruzsa-plun}]\label{plun-ruz} Let $A$ be contained in an abelian group, and suppose that $A$ has doubling $K$. Then $|kA - lA| \leq K^{k+l}|A|$ for all $k,l \geq 0$.\end{proposition}

\section{Finding a good model.}
\label{sec2} \ni Our aim in this section is to prove Proposition \ref{prop2}. If $A$ is a subset of some abelian group then by an $s$-\textit{model} for $A$ we mean a pair $(A',G')$, where $G'$ is an abelian group and $A' \subseteq G'$ is Freiman $s$-isomorphic to $A'$. If $G'$ is finite then we call $(A',G')$ a \textit{finite model}, and by the \textit{size} of such a model we mean simply the cardinality of $G'$. A \textit{minimal} $s$-model is a model for $A$ with the minimal size. The existence of a minimal $s$-model is a trivial consequence of our first lemma.

\begin{lemma}[Finite sets have finite models]\label{ff}
Let $A$ be a finite subset of some abelian group, and let $s \geq 2$ be an integer. Then $A$ has a finite $s$-model.
\end{lemma}
\proof Suppose that $G$ is the group generated by $A$. By the structure theorem for finitely-generated abelian groups, $G$ is isomorphic to $H \times \Z^k$ for some non-negative integer $k$. Regard $A$ as a subset of $H \times \Z^k$. Then there is certainly some integer $N$ such that $A \subseteq H \times [-N,N]^k$. The projection map 
\[ \pi : H \times \mathbb{Z}^k \rightarrow H \times (\Z/4sN\Z)^k\]
induces a Freiman $s$-isomorphism on $A$, and therefore $A$ has an $s$-model inside the finite group $H \times (\Z/4sN\Z)^k$.\endproof\vs

\ni Our next lemma is a simple consequence of \cite[Lemma 11]{green-ruzsa}.
\begin{lemma}\label{lem3A} Let $\epsilon \in (0,\frac{1}{200})$. Let $A \subseteq G$ be a set with doubling constant $K$. Write $D = A - A$, and suppose that
\[ \E 1_A  \leq  K^{-2}\epsilon^{4K^2}.\] Then there is a non-trivial character $\gamma \in G^{\ast}$ such that $|\widehat{1}_D(\gamma)| \geq (1 - \epsilon)\E 1_D$.
\end{lemma}
\proof Apply \cite[Lemma 11]{green-ruzsa}. This tells us that there is $\gamma \neq \gamma_0$ such that $|\widehat{1}_D(\gamma)| \geq (1 - \eta)\E 1_D$, where 
\[ \eta  =  9K^{-2} (\E 1_D)^{1/2K^2} \log (1 / \E 1_D).\] By a result of the second author \cite{ruzsa-sumdiff} we have $\E 1_D \leq K^2 \E 1_A$ (this is also a special case of Proposition \ref{plun-ruz}). It is easy to check that the condition on $\E 1_A$ in the statement of the lemma is enough to ensure that $\eta \leq \epsilon$.\endproof\\[11pt]
The next lemma is a slight variation on \cite[Lemma 7]{green-ruzsa}. 
\begin{lemma}\label{lem4} Let $A \subseteq G$, and let $\psi : G \rightarrow \mathbb{Z}/q\mathbb{Z}$ be a homomorphism of groups. Suppose that 
\[ |(A - A) \setminus \psi^{-1}[b, b+l]| \; < \; |A|/2\] for some $b \in \mathbb{Z}/q\mathbb{Z}$ and some $l < q/3$. Then $\psi(A) \subseteq [x,x+l]$ for some $x$.
\end{lemma}
\proof By translating if necessary we may assume that $0 \in A$, in which case $A \subseteq A - A$, and that the longest gap in $\psi(A)$ is the interval $[k,-1]$. Write $A' = A \cap \psi^{-1}[b,b+l]$. Observing that $|A \setminus \psi^{-1}[b,b+l]| < |A|/2$, we see that $|A'| > |A|/2$. Thus for an arbitrary $a \in A$ the set $A' - a$, which is a subset of $A - A$, cannot be contained in $(A - A) \setminus \psi^{-1}[b,b+l]$. This means that there is at least one $a' \in A'$ for which $\psi(a' - a) \in [b,b+l]$, from which it follows that $\psi(a) \in \psi(a') - [b,b+l] \subseteq [-l,l]$.\vs

\ni We have proved, then, that $\psi(A) \subseteq [-l,l]$, which means that there is a gap in $\psi(A)$ of length $(q - l) - l > l$. However, we assumed that the longest gap in $\psi(A)$ was of the form $[k,-1]$, from which we may conclude that $\psi(A) \cap [-l,0)$ is empty. Thus $\psi(A) \subseteq [0,l]$, which is what we wanted to prove (remember we started by subjecting $A$ to a translation).\endproof
\begin{lemma}\label{lem5} Let $S \subseteq G$ be any set. Suppose that $\kappa,\delta$ are real numbers such that $0 < \kappa < 1$, $0 < \delta < 1/2$. Suppose that $\gamma \neq \gamma_0$ is a character for which $|\widehat{1}_S(\gamma)| \geq (1 - 4\kappa\delta^2)\E 1_S$. Write $\mbox{\emph{ord}}(\gamma) = q$, and write $\psi : G \rightarrow \mathbb{Z}/q\mathbb{Z}$ for the homomorphism obtained by composing $\gamma$ with the map which sends $e^{2\pi ij/q}$ to $j\pmod{q}$. Then there are $b$ and $l$, $l < \delta q$, such that 
\[ |S \setminus \psi^{-1}[b,b+l]| \; \leq \; \kappa |S|.\]
\end{lemma}
\proof Write $\widehat{1}_S(\gamma) = e^{2\pi i \lambda/q}(1 - \eta)\E 1_S$ with $\lambda \in \R$ and $0 \leq \eta \leq 4\kappa \delta^2$. Then
\[ (1 - \eta)|S| \; = \; \sum_{s \in S} e^{2\pi i(\psi(s) - \lambda)/q} \; = \; \sum_{s \in S} \cos \frac{2\pi (\psi(s) - \lambda)}{q}.\]
Each summand is at most $1$ and, moreover, for those $s$ such that $\psi(s) \notin (\lambda - \delta q/2,\lambda + \delta q/2)$ it is at most $\cos \delta \pi < 1 - 4\delta^2$. It follows that the number of such $s$ is less than $\eta |S|/4\delta^2 \leq \kappa |S|$, as required.\endproof\vs

\ni The next result links the above three lemmas.
\begin{lemma}\label{lem6}
Suppose that $A \subseteq G$ has doubling constant $K$. Let $\delta \in (0,\frac{1}{20})$, and suppose that $\E 1_A \leq (\delta/K)^{10K^2}$. Then there is $q \geq 2$ and a homomorphism $\psi : G \rightarrow \mathbb{Z}/q\mathbb{Z}$ such that $\psi(A) \subseteq [x, x+ l]$ for some $x$ and some $l < \delta q$.
\end{lemma}
\ni The conditions of Lemma \ref{lem3A} are satisfied with $\epsilon = \delta^2/K^2$. This shows that if $D = A - A$ then there is $\gamma \neq \gamma_0$ such that $|\widehat{1}_D(\gamma)| \geq (1 - \epsilon)\E 1_D$. This means that the hypotheses of Lemma \ref{lem5} are satisfied with $\kappa = 1/4K^2$, which means that for some $b$ and $l < \delta q$ we have
\[ |D \setminus \psi^{-1}[b,b+l]| \leq  \kappa |S|  =  \frac{1}{4K^2}|A - A|  <  |A|/2,\] this last step following by another application of \cite{ruzsa-sumdiff}. Finally, we may apply Lemma \ref{lem4} to get the desired conclusion.\endproof\vs

\ni We may now supply a proof of Proposition \ref{prop2}. Recall that $A \subseteq G$ is a set with doubling constant $K$, and that $s \geq 2$ is an integer.\vs

\ni\textit{Proof of Proposition \ref{prop2}.} We may suppose that $(A,G)$ is already a minimal $s$-model; the existence of such a model follows from Lemma \ref{ff}. Suppose for a contradiction that $\E 1_A < (10sK)^{-10K^2}$.\vs

\ni By Lemma \ref{lem6} there is $q \geq 2$ and a homomorphism $\psi : G \rightarrow \mathbb{Z}/q\mathbb{Z}$ such that $\psi(A) \subseteq [x,x+l]$ for some $x$ and some $l < q/4s$. By translating $A$ if necessary, we may assume that $x = 0$. Now let $H = \ker \psi$ and let $z \in \psi^{-1}(1)$, so that $H$ and $z$ together generate $G$. Let $G' = H \times \mathbb{Z}/(q-1)\mathbb{Z}$, and consider the map $\theta : G \rightarrow G'$ defined as follows. If $g \in G$, write $g = h + \lambda z$ with $h \in H$ and $0 \leq \lambda < q$. Then $\theta(g) = (h,\lambda\pmod{q-1})$. $\theta$ is certainly not a group homomorphism, but it does induce a Freiman $s$-isomorphism on $A$. To see this, suppose that $a_1,\dots,a_s,a'_1,\dots,a'_s$ are $2s$ elements of $A$ satisfying
\[ a_1 + \dots + a_s \; = \; a'_1 + \dots + a'_s.\]
Write $a_i = h_i + \lambda_i z$, $a'_i = h'_i + \lambda'_i z$ where $0 \leq \lambda_i \leq q/2s$. Then we must actually have 
\[ \lambda_1 + \dots + \lambda_s \; = \; \lambda'_1 + \dots + \lambda'_s,\] and so 
\[ \theta(a_1) + \dots + \theta(a_s) \; = \; \theta(a'_1) + \dots + \theta(a'_s).\]
The proof that $\theta^{-1}$ is an $s$-homomorphism is very similar. \vs

\ni We have shown that $A$ has a model in $G'$, which is contrary to our assumption that $(A,G)$ was a minimal model.\endproof

\section{The Bogolyubov-Chang argument.}

\ni Let $A \subseteq G$ be a set with doubling $K$. Suppose that it is a minimal $8$-model. Then it follows from Proposition \ref{prop2}, as proved in the preceding section, that $|G| \leq (80K)^{10K^2}|A|$.\vs

\ni In this section we will show that $2A - 2A$ contains a large Bohr set in $G$. In \S \ref{sec4} we will prove a result on the structure of Bohr sets, and then use this to get information when $(A,G)$ is not a minimal model.\vs

\begin{proposition}\label{changbog}
Let $A \subseteq G$, and suppose that $A$ has doubling $K$. Then $2A - 2A$ contains some Bohr neighbourhood $B(\Gamma,\delta)$, 
where $|\Gamma| \leq 8K\log(1/\E 1_A)$ and $\delta \geq 
\left(48K\log(1/\E 1_A)\right)^{-1}$.
\end{proposition}

\ni The proof of this result, when $G = \mathbb{Z}/N\mathbb{Z}$, is a combination of ideas of Bogolyubov \cite{bogolyubov} and Chang \cite{chang}. Chang built upon earlier ideas of Rudin \cite{rudin}. Generalising to the case of a general abelian group $G$ is completely straightforward. However, we would like to take the opportunity to record an alternative approach to the Chang-Rudin part of the argument, the heart of which is Proposition \ref{chang} below. Doing so has the advantage of keeping the paper relatively self-contained.\vs

\ni We begin with some notation and definitions. Let $G$ be an abelian group, and let $\Phi = \{\phi_1,\dots,\phi_d\} \subseteq \widehat G$ be a set of characters. Write $\langle \Phi \rangle$ for the \textit{cube} spanned by $\Phi$, that is the collection of all characters $\phi_1^{\epsilon_1}\dots \phi_d^{\epsilon_d}$ where $\epsilon_j \in \{-1,0,1\}$ for all $j$. We call $d$ the \textit{dimension} of such a cube. The set $\Phi$ is said to be
\textit{dissociated} if the only solution to $\phi_1^{\epsilon_1}\dots \phi_d^{\epsilon_d} = 1$ with $\epsilon_j \in \{-1,0,1\}$ is 
the trivial one, in which $\epsilon_j = 0$ for all $j$.
\begin{proposition}[Chang]\label{chang}  Let $\rho,\alpha \in [0,1]$, let $A 
\subseteq G$ be a set of size $\alpha |G|$ and let $\Gamma \subseteq \widehat G$ be the set of all $\gamma$ for which $|\widehat{1}_A(\gamma)| \geq \rho \E1_A$. Let 
$\Phi$ be a dissociated subset of $\Gamma$. Then $|\Phi| \leq 2\rho^{-
2}\log(1/\E 1_A)$.
\end{proposition}
\ni Throughout what follows we will suppose that $\Phi = \{\phi_1,\dots,\phi_d\}$ is a dissociated set of characters. We will consider functions of the form
\begin{equation}\label{touse} f(x) \; = \; \sum_{j=1}^d c_j \re(\omega_j \phi_j(x)),\end{equation} where $c_j \in \mathbb{R}$ and $|\omega_j| = 1$ for $j = 
1,\dots,d$. We shall show, by remodelling a classical 
technique of Bernstein which is nearly 80 years old, that such functions behave rather like
sums of independent random variables. We will then derive Proposition \ref{chang} from this observation.
\begin{lemma}\label{lem3} Let $f$ be given by \eqref{touse}. Then we have $\E f^2  = \frac{1}{2}\sum_j c_j^2$.
\end{lemma}
\proof If $i \neq j$ then one may check that
\[ \E_{x \in G} \re(\omega_i \phi_i(x))\re(\omega_j \phi_j(x)) \; = \; 0.\] This follows using the orthogonality relations $\E_{x \in G} \phi_i(x)\phi_j^{-1}(x) = 0$ and $\E_{x \in G} \phi_i(x)\phi_j(x) = 0$, this second relation being a consequence of the fact (a very minor consequence of dissociativity) that $\phi_j \neq \phi_i^{-1}$. It is just as easy to check that 
\[ \E_{x \in G} \re(\omega_i \phi_i(x))\re(\omega_i \phi_i(x))  =  1/2.\]
The lemma follows by combining these two pieces of information in the obvious way.\endproof
\begin{proposition}\label{prop1}
Let $t \in \mathbb{R}$ and let $f$ be given by \eqref{touse}. Then we have the inequality
\[  \E e^{t f(x)}  \leq  e^{t^2 \E f^2}.\]
\end{proposition}
\proof For any real $y$ satisfying $|y| \leq 1$ one has the inequality $e^{ty} \leq \cosh(t) + y \sinh(t)$. Thus
\begin{equation}\label{starr} \E_{x \in G} e^{t f(x)}  \leq  \E_{x \in G} 
\prod_{j = 1}^k \left( \cosh(t c_j) + \sinh(t c_j)\re(\omega_j\phi_j(x)\right).\end{equation}
Write $\re(\omega_j\phi_j(x)) = \frac{1}{2}(\omega_j\phi_j(x) + \overline{\omega_j}\phi_j^{-1}(x))$, and multiply out to get a linear combination of terms $\phi_1^{\epsilon_1}(x)\dots\phi_k^{\epsilon_k}(x)$, $\epsilon_i \in \{-1,0,1\}$. Since $\Phi$ is dissociated, only the term with $\epsilon_1 = \dots = \epsilon_k = 0$ does not vanish when we take expectations over $x$. The coefficient of this term in \eqref{starr} is $\prod_{j = 1}^k \cosh(t c_j)$. Therefore, using Lemma \ref{lem3} and the elementary inequality $\cosh(u) \leq e^{u^2/2}$ we obtain
\[ \E e^{t f(x)}  \leq  \prod_{j = 1}^k \cosh(t 
c_j) \leq  \exp\big(\textstyle\frac{1}{2}\displaystyle t^2\sum_{j=1}^d c_j^2 \big) = 
e^{t^2 \E f^2},\] as desired.\endproof\\[11pt]
\textit{Proof of Proposition \ref{chang}.} Set \[ f(x)  = \re \big(\sum_{j=1}^d \widehat{1}_A(\phi_j) \phi_j(x)\big).\]
This is certainly a function of the form \eqref{touse}, in 
which $c_j = |\widehat{1}_A(\phi_j)|$.
Observe also that $\widehat{f}(\gamma) = \widehat{1}_A(\gamma)/2$ if $\gamma \in \Phi \cup -\Phi$ 
and $\widehat{f}(\gamma) = 0$ otherwise. This implies that
\begin{equation}\label{useinasec}(\E 1_A)\E_{x \in A} f(x) =  \E_{x \in G} f(x)1_A(x)  = \sum_{\gamma \in \widehat{G}} 
\widehat{f}(\gamma)\overline{\widehat{1}_A(\gamma)}  = 2\sum_{\gamma \in \widehat{G}} |\widehat{f}(\gamma)|^2  =  
2\E f^2.\end{equation}
Proposition \ref{prop1} and 
\eqref{useinasec} give
\[
\frac{1}{\E 1_A}e^{t^2 \E f^2}  \geq   
\frac{1}{\E 1_A} \E_{x \in G} e^{t f(x)} \geq 
 \E_{x \in A} e^{tf(x)} \geq  
\exp\big(t\E_{x \in A} f(x)\big) =  
\exp\big(\frac{2t \E f^2}{\E 1_A}\big),\]
the third inequality being a consequence of the convexity property
\[   \frac{e^{x_1} + \dots + e^{x_n}}{n} \geq \exp \bigg(\frac{x_1 + \dots + x_n}{n}\bigg)\] (the weighted A.M.-- G.M. inequality).
Choosing $t = 1/\E1_A$ gives \begin{equation}\label{eq11} 1/\E1_A  \geq  
\exp\big( \E f^2 / (\E1_A)^2 \big).\end{equation} Now using Lemma \ref{lem3} and our assumption that $\Phi$ is a subset of $\Gamma$, the set of all $\gamma \in G^{\ast}$ such that $|\widehat{1}_A(\gamma)| \geq \rho \E 1_A$, we have that 
\[ 2\E f^2  =  \sum_j c_j^2  =  
 \sum_{\phi \in \Phi}|\widehat{1}_A(\phi)|^2  \geq  
d\rho^2 (\E 1_A)^2.\] Comparing with \eqref{eq11} concludes the proof of Proposition \ref{chang}.\endproof\\[11pt]
\noindent\textit{Proof of Proposition \ref{changbog}.}
Recall that $A \subseteq G$ has doubling $K$. Since the convolution $1_A \ast 1_A$ is supported on $A + A$, the Cauchy-Schwarz inequality gives
\begin{equation}\label{L4A} \sum_{\gamma \in \widehat{G}} |\widehat{1}_A(\gamma)|^4  = 
\E_{x \in G} (1_A \ast 1_A)(x)^2 1_{A + A}(x)  \geq  \frac{\big(\E_{x \in G} 1_A \ast 1_A(x)\big)^2}{\E 1_{A + A}} =  \frac{(\E 1_A)^4}{\E 1_{A + A}}  = 
\frac{(\E 1_A)^3}{K}.\end{equation}
Now let $\Gamma$ be the set of all $\gamma$ for which $|\widehat{1}_A(\gamma)| \geq 
\E1_A/2\sqrt{K}$. We claim that $2A - 2A$ contains $B(\Gamma,\frac{1}{6})$. Indeed if 
$x \in B(\Gamma,\frac{1}{6})$ then $\re\gamma(x) \geq 1/2$ for all $\gamma \in \Gamma$. It follows from the inversion formula that
\begin{eqnarray*}
(1_A \ast 1_A \ast 1_A \ast 1_A)(x) & = & \sum_{\gamma \in \widehat G} |\widehat{1}_A(\gamma)|^4 \overline{\gamma(x)} \\ & = & \sum_{\gamma \in \widehat G} |\widehat{1}_A(\gamma)|^4 \re \gamma(x) \\ & > & \frac{1}{2}\sum_{\gamma \in \Gamma} 
|\widehat{1}_A(\gamma)|^4 - 2 \sum_{\gamma \notin \Gamma}|\widehat{1}_A(\gamma)|^4\\ & \geq &  
\frac{(\E 1_A)^3}{2K} - 2 \sup_{\gamma \notin \Gamma}|\widehat{1}_A(\gamma)|^2 \sum_{\gamma} 
|\widehat{1}_A(\gamma)|^2 \\ & \geq & 0.
\end{eqnarray*}
(Note that the inequality on the third line is strict since $\re \gamma(x) = 1$ when $\gamma = \gamma_0$ is the trivial character.)
It follows that $B(\Gamma,\frac{1}{6})$ is indeed contained in $2A - 2A$. Parseval's identity tells us that $|\Gamma| 
\leq 4K/\E 1_A$, but Proposition \ref{changbog} claims that $2A - 2A$ in fact contains a Bohr set of much smaller dimension than this. To obtain such a result, we apply Proposition \ref{chang}.
Let $\Phi = \{\phi_1,\dots,\phi_d\}$ be the largest dissociated subset of $\Gamma$. By Proposition \ref{chang} we have $d \leq 8K\log(1/\E 1_A)$. Now if $\gamma \in \Gamma \setminus \Phi$ then there must be a relation of the form $\gamma\phi_1^{\epsilon_1}\dots\phi_d^{\epsilon_d} = 0$, since otherwise $\Phi \cup \gamma$ would be dissociated, contrary to the assumption that $\Phi$ is the \textit{maximal} dissociated subset of $\Gamma$. Thus $\Gamma$ is contained in $\langle \Phi \rangle$, the cube spanned by $\Phi$.\vs

\ni We claim that $B(\Gamma,\frac{1}{6})$ contains 
$B(\Phi,1/6d)$, which immediately implies Proposition \ref{changbog}. Indeed, any $\gamma \in \Gamma$ can be 
written as $\phi_1^{\epsilon_1}\dots\phi_d^{\epsilon_d}$ with $\epsilon_j \in 
\{-1,0,1\}$. Supposing then that $x \in B(\Phi,1/6d)$, so that $(2\pi)^{-1}|\arg(\phi_i(x))| \leq 1/6d$ for all $i$, we have
\[ |\arg(\gamma(x))| \leq \sum_{i = 1}^d |\arg(\phi_i(x))| \leq 2\pi \cdot \frac{1}{6}.\]
This confirms the claim, and hence completes the proof of Proposition \ref{changbog}.\endproof

\section{The structure of Bohr sets.}\label{sec4}

\ni In this section we shall apply results from the geometry of numbers to elucidate the structure of Bohr sets.\vs

\ni We begin by recalling \textit{Minkowski's second theorem}, a proof of which may be found in \cite[Ch. VIII, Thm. V]{cassels}. If $\Lambda \subseteq \R^d$ is a lattice and if $Q$ is a centrally symmetric closed convex body, we define the \textit{successive minima} $\lambda_1,\dots,\lambda_d$ of $Q$ with respect to $\Lambda$ by
\[ \lambda_j := \inf \{ \lambda : \lambda Q \cap \Lambda \;\; \mbox{contains $j$ linearly independent vectors}\}.\]

\begin{proposition}[Minkowski's second theorem] We have the bound
\[ \lambda_1 \lambda_2 \dots \lambda_d |Q| \leq 2^d \det(\Lambda).\]
\end{proposition}

\begin{proposition}\label{prop4.2}
Let $\Gamma \subseteq \widehat{G}$, $\Gamma = \{\gamma_1,\dots,\gamma_d\}$, be a set of $d$ characters, and let $\rho \in (0,1/4)$. Then the Bohr set $B(\Gamma,\rho)$ contains a proper coset progression $P + H$ with dimension $d$ and size at least $(\rho/d)^d|G|$. 
\end{proposition}
\proof Define \[ H := \bigcap_{j=1}^d \ker \gamma_j.\] If $x = (x_1,\dots,x_d)$ is an element of $\T^d$, written with $|x_j| \leq 1/2$ for all $j$, then we write \[ \Vert x \Vert_{\infty} := \sup_{1 \leq j \leq d} |x_j|.\] Consider the map $\phi : G \rightarrow \T^d$ defined by
\[ \phi(x) = (\arg \gamma_1(x),\dots, \arg \gamma_d(x)).\]
The image $\phi(G)$ is a subgroup of $\T^d$, and the Bohr set $B(\Gamma,\rho)$ is simply the inverse image under $\phi$ of the cube $\rho Q := \{x : \Vert x \Vert_{\infty} \leq \rho\}$.
Let $\Lambda \subseteq \R^d$ be the subgroup $\phi(G) + \Z^d$; observe that  
\begin{equation}\label{det-lambda-formula} \det(\Lambda) = |H|/|G|.\end{equation}
With a slight abuse of notation, we have
\[ \phi(G) \cap \rho Q = \Lambda \cap \rho Q.\] On the right hand side, $Q$ now refers to the cube
\[  Q := \{ x \in \R^d : \Vert x \Vert_{\infty} \leq 1\}.\]
This being a centrally symmetric convex body, we may apply Minkowski's second theorem. Writing $\lambda_1,\dots,\lambda_d$ for the successive minima of $Q$ with respect to $\Lambda$, this theorem and \eqref{det-lambda-formula} tell us that 
\begin{equation}\label{mink-conseq}
\lambda_1 \lambda_2 \dots \lambda_d \leq |H|/|G|.
\end{equation}
Now we may choose, inductively, linearly independent vectors $b_1,\dots,b_d \in \Lambda$ such that $b_j \in \lambda_j Q$, which means that $\Vert b_j \Vert_{\infty} \leq \lambda_j$. It is clear that all of the vectors
\[ \{ l_1 b_1 + \dots + l_d b_d : -L_j \leq l_j \leq L_j, l_j \in \Z\}\] lie in $\rho Q$, where $L_j := \lfloor \rho/d\lambda_j\rfloor$. Let $v_j \in G$ be arbitrary elements for which $\phi(v_j) \equiv b_j \md{\Z^d}$, and define
\[ P := \{ l_1 v_1 + \dots + l_d v_d : -L_j \leq l_j \leq L_j\}.\]
Then $P + H \subseteq B(\Gamma,\rho)$, and it is easy to see that $P + H$ is a proper coset progression. It remains to give a lower bound for $|P + H|$. To this end, note that the number of integers $l_j$ with $-L_j \leq l_j \leq L_j$ is at least $\rho/d\lambda_j$, and so by \eqref{mink-conseq} we have
\[ |P| \geq \big(\frac{\rho}{d}\big)^d \frac{1}{\lambda_1 \dots \lambda_d} \geq \big(\frac{\rho}{d}\big)^d \frac{|G|}{|H|}.\]
The proposition follows immediately.\endproof
\section{A covering argument.}

\ni In this section we conclude the proof of Theorem \ref{mainthm} by employing a covering argument of Chang \cite{chang}. Let us begin by combining the results of the last few sections.\vs

\ni Let $A \subseteq G$ be a set with doubling $K$. Suppose that $(A',G')$ is a minimal $8$-model, and that $\pi : A \rightarrow A'$ is a Freiman $8$-isomorphism. Then it follows from Proposition \ref{prop2}, which was proved in \S \ref{sec2}, that $|G'| \leq (80K)^{10K^2}|A|$. Applying Proposition \ref{changbog} to $A'$ and noting that $\E 1_{A'} \geq (80K)^{-10K^2}$, we see that $2A' - 2A'$ contains a Bohr set $B(\Gamma,\rho)$ with $|\Gamma| \leq 2^9 K^3 \log(K+2)$ and $\rho \geq (2^{12} K^3 \log(K+2))^{-1}$. By Proposition \ref{prop4.2}, this Bohr set in turn contains a proper coset progression $P + H$ with dimension at most $2^9 K^3 \log(K+2)$ and size satisfying
\begin{equation}\label{size} |P + H| \geq \exp(-2^{14} K^3 \log^2 (K+2)) |G'| \geq \exp( -2^{14} K^3 \log^2(K+2)) |A|.\end{equation}
Now coset progressions, together with their size and dimension, are preserved under Freiman $2$-isomorphisms. Since the Freiman $8$-isomorphism $\pi^{-1}$ on $A'$ induces a $2$-isomorphism on $2A' - 2A'$, it follows that $2A - 2A$ contains a coset progression, which we will also call $P + H$, with dimension at most $2^9 K^3 \log (K+2)$ and size satisfying the lower bound \eqref{size}.\vs

\ni At this point the model $(A',G')$, which we considered so that we could do harmonic analysis, has served its purpose. Henceforth we will work only with $A$ and the coset progression $P + H$ contained in $2A - 2A$. The next proposition deals with this situation, and Theorem \ref{mainthm} is a straightforward consequence of it.

\begin{proposition}[Chang] Suppose that $A$ is a subset of an abelian group with doubling $K$ and that $2A - 2A$ contains a proper coset progression $P + H$ 
of size $\eta |A|$ and dimension $d$. Then $A$ is contained in a coset progression of 
size at most $2^d(K^4\eta^{-1})^{5K}|A|$ and dimension at most $d + 
4K\log(K^4/\eta)$.
\end{proposition}
\proof We describe an algorithm for selecting some non-negative integer $t$ and 
subsets $S_i$, $i \leq t$, of $A$. Set $P_0 := P + H$. Let $R_0$ be a maximal subset 
of $A$ for which the translates $P_0 + x$, $x \in R_0$, are all disjoint. If 
$|R_0| \leq 2K$ then set $t = 0$ and $S_0 = R_0$, and terminate the algorithm. 
Otherwise take $S_0$ to be any subset of $R_0$ of cardinality $2K$, and set $P_1 
= P_0 + S_0$. Take $R_1$ to be a maximal subset of $A$ for which the translates 
$P_1 + x$, $x \in R_1$, are all distinct. If $|R_1| \leq 2K$ then set $t = 1$ 
and $S_1 = R_1$ and terminate the algorithm. Otherwise choose $S_1 \subseteq 
R_1$ with $|S_1| = 2K$ and set $P_2 = P_1 + S_1$. Continue in this way.\vs

\ni We claim that this is a finite algorithm, and that in fact $t \leq 
\log(K^4/\eta)$. Indeed the fact that the translates $P_i + x$, $x \in S_i$, are 
all disjoint means that $|P_{i+1}| = |P_i||S_i|$ for $i \leq t - 1$. It follows 
that
\begin{equation}\label{eqq7} |P_t| \; \geq \; |P + H||S_0|\dots|S_{t-1}|  \geq  
\eta(2K)^t|A|.\end{equation}
Observe, however, that
\[ P_t  \subseteq  P + H + A + A + \dots + A,\] where there are $t$ copies of 
$A$. Since $P + H \subseteq 2A - 2A$ this means that $P_t \subseteq (t+2)A - 2A$, 
and hence by Proposition \ref{plun-ruz} we have $|P_t| \leq K^{t + 4}|A|$. 
Comparison with (\ref{eqq7}) proves the claim.\\[11pt]
Let us examine what happens when the algorithm finishes. Then we have a set $R_t 
\subseteq A$, $|R_t| \leq 2K$, which is maximal subject to the translates $P_t + 
x$, $x \in R_t$, being disjoint. In other words if $a \in A$ then there is $x 
\in R_t$ such that $(P_t + a) \cap (P_t + x) \neq \emptyset$, and so
\begin{equation}\label{eq47} A  \subseteq P_t - P_t + R_t  \subseteq 
(P - P) + (S_0 - S_0) + \dots + (S_{t-1} - S_{t-1}) + R_t + H.\end{equation}
Now if $S = \{s_1,\dots,s_m\}$ is a subset of an abelian group then the cube 
\[ \overline{S} := \{ \epsilon_1 s_1 + \dots + \epsilon_m s_m : \epsilon_j \in \{-1,0,1\}\}\] 
is a multidimensional progression of dimension $|S|$ and size at most $3^{|S|}$, 
and it contains the set $S - S$. It follows from (\ref{eq47}) that $A \subseteq 
Q + H$, where $Q$ is the progression
\[ Q  =  P - P + \overline{S}_0 + \dots + \overline{S}_{t-1} + 
\overline{R}_t.\]
The dimension of $Q$ satisfies
\begin{eqnarray*} \mbox{dim}(Q) & \leq & \mbox{dim}(P) + \sum_{i = 0}^{t-1} 
|S_i| + |R_t| \\ & \leq & d + 2K(t+1) \\ & \leq & d + 
4K\log(K^4/\eta).\end{eqnarray*}
To estimate the size of $Q + H$, note that the properness of $P$ implies that $|P - 
P| = 2^d|P|$. Hence
\begin{eqnarray*} |Q + H| & \leq & |H||P - P|\cdot \prod_{i = 0}^{t-1} 3^{|S_i|} \cdot 
3^{|R_t|} \\ & \leq & 2^d3^{2K(t+1)}|P||H| \\ & \leq & 
2^d3^{4K\log(K^4/\eta)}K^4|A|\\ & \leq & 
2^d\left(\frac{K^4}{\eta}\right)^{5K}|A|,\end{eqnarray*}
the penultimate step following from Proposition \ref{plun-ruz} and the fact that 
$P + H \subseteq 2A - 2A$.\endproof\vs

\ni By the remarks at the start of the section we may apply this proposition with $d \ll K^3 \log (K+2)$ and $\eta \geq \exp(-C K^3 \log^2 (K+2))$. This leads immediately to Theorem \ref{mainthm}.\endproof

\section{Further remarks on models.}
\ni We feel that the notion of a model is one that could be investigated further. In this section we collect a few further remarks on $s$-models, restricting ourselves for simplicity to the case $s = 2$. When we talk about isomorphisms or models, we mean $2$-isomorphisms and $2$-models. We will also, in this section, suppose that the doubling constant $K$ is larger than $1000$ so as to avoid having to make tedious estimates valid for all $K$, and that $|A| > n_0(K)$.\vs

\ni Proposition \ref{prop2} tells us that if $A$ is a subset of an abelian group with doubling $K$, then $A$ has a model of size at most $e^{20 K^2 \log K}|A|$. Our first two results show that if $A$ is either a subset of $\mathbb{F}_2^{m}$ or a set of integers then a much smaller model can be found.
\begin{proposition}\label{lem6.1} Suppose that $A \subseteq \mathbb{F}_2^m$ has doubling $K$. Then $A$ has a model of size at most $K^4 |A|$.
\end{proposition}
\proof Suppose that $m$ is the smallest positive integer such that $A$ is isomorphic to a subset of $\mathbb{F}_2^m$. Now by Proposition \ref{plun-ruz} we have $|2A - 2A| \leq K^4|A|$. Thus if $|2^m| > K^4|A|$ then there is some $x \in \mathbb{F}_2^m \setminus (2A - 2A)$. Let $\phi : \mathbb{F}_2^{m} \rightarrow \mathbb{F}_2^{m-1}$ be any linear map with kernel $\{0,x\}$. Clearly $\phi$ induces a Freiman homomorphism on $A$. If, however, we have
\[ \phi(a_1) + \phi(a_2) \; = \; \phi(a'_1) + \phi(a'_2)\] then
$a_1 + a_2 - a'_1 - a'_2 \in \ker \phi = \{0,x\}$. This means, since $x \notin 2A - 2A$, that $a_1 + a_2 = a'_1 + a'_2$, which implies that $\phi$ in fact induces a Freiman isomorphism on $A$. This is contrary to the supposed minimality of $m$.\endproof
\begin{proposition}\label{lem6.2} Suppose that $A \subseteq \mathbb{Z}$ has doubling $K$. Then $A$ has a model of size at most $100K^6\log K |A|$.
\end{proposition}
\proof Write $n := |A|$. By \cite[Theorem 3]{green-ruzsa}, $A$ is isomorphic to a subset of the interval $[1,\frac{1}{2}f(K)]$, where $f(K) = e^{8K^2\log K}$. Suppose for a contradiction that none of the projections $\pi_m : \mathbb{Z} \rightarrow \mathbb{Z}/m\mathbb{Z}$, $m \leq 100 K^6 \log K \cdot n$, induces a Freiman isomorphism of $A$. Then, by much the same reasoning we used in Proposition \ref{lem6.1}, the set $2A - 2A$ must contain a multiple $\lambda_m m$ for all $m \leq 100K^6\log K \cdot n$. Set $L = 50K^6\log K$ , $X = 2f(K)/L$ and let
\[ S  =  \left\{ m \in [Ln,2Ln] : p | m \Rightarrow p \geq X.\right\}.\]
We claim that the elements $\lambda_s s$, $s \in S$, are all distinct. Indeed if $s,s' \in S$ are distinct then \[ \mbox{lcm}(s,s') \; \geq \; X \min(s,s') \; \geq \; XLn \; > \; f(K)n,\]
whereas $2A - 2A \subseteq [-f(K)n,f(K)n]$. However (for large $n$) we have the estimate
\[ |S|  >  \frac{1}{2}L \prod_{p \leq X}\left(1 - \frac{1}{p}\right) \; \geq \; \frac{Ln}{4\log X} \; > \; K^4 n.\]
Thus $|2A -2A| > K^4 n$, which is contrary to the estimate $|2A - 2A| \leq K^4 n$ furnished by Proposition \ref{plun-ruz}.\endproof\vs

\ni Observe that the model was constructed in $\mathbb{Z}/m\mathbb{Z}$, where $m$ has no small prime factors. It would be interesting to know whether, in fact, any set of $n$ integers with doubling $K$ has a model in $\mathbb{Z}/p\mathbb{Z}$ with $p$ a smallish prime, say $p \leq K^{100}n$. We cannot decide whether or not this is so, though it is natural to suppose once again that $A \subseteq [1,\frac{1}{2}f(K)n]$ and to consider the projections $\pi_p : \mathbb{Z} \rightarrow \mathbb{Z}/p\mathbb{Z}$, $p \in [Ln,2Ln]$. If none of these induces a Freiman homomorphism then $2A - 2A$ must contain a multiple $\lambda_p p$, where $\lambda_p \leq f(K)/L$. Perhaps, for a suitable choice of $L$, this is at odds with the fact that $|4A - 4A| \leq K^8 n$ (another consequence of Proposition \ref{plun-ruz}). This line of thinking leads us to formulate the following question.
\begin{question}
Let $P$ and $X$ be positive constants. Suppose that $S \subseteq \mathbb{Z}$ is a set such that for each prime $P \leq p < 2P$ there is an element $s_p \in S$ of the form $\lambda_p p$, where $\lambda_p \leq X$ is a positive integer. Find a lower bound for $|S + S|$.
\end{question}
\ni Perhaps it is true that $|S + S| \gg P/(\log X)^{\alpha}$ for some absolute $\alpha$, which would imply the existence of small models modulo a prime as outlined above.\\[11pt]
Propositions \ref{lem6.1} and \ref{lem6.2} might lead one to believe that \textit{any} set $A$ with doubling $K$ has a model of size $K^{100} |A|$. We conclude this section by showing that this is not so.
\begin{proposition}
Fix a constant $K$. Then there are infinitely many $n$ with the following property. There exists a set $A$ \emph{(}in some abelian group $G$\emph{)} with cardinality $n$ and doubling at most $K$ but which has no model in any group $G'$ satisfying $|G'| \leq e^{\frac{1}{6}\sqrt{K \log K}}n$.
\end{proposition}
\proof Set $X = \lfloor \frac{1}{3}\sqrt{K \log K}\rfloor$, and let $p_1 = 2,p_2 = 3,\dots,p_k$ be the primes less than $X$. Let $Q$ be any prime larger than $X$, and write $G = \mathbb{Z}/Q\mathbb{Z} \times \prod_{i = 1}^k \mathbb{Z}/p_i\mathbb{Z}$. Consider the set $A$ consisting of all $(k+1)$-tuples $(x,x_1,\dots,x_k)$ in which at most one of the $x_i$ is non-zero. We claim that any group $G'$ which contains a model for $A$ must contain elements of order $p_1,\dots,p_k$ and $Q$. Indeed, suppose that $A' \subseteq G'$ and that $\psi : A \rightarrow A'$ is a Freiman isomorphism.
Let $t = (0,0,\dots,0,1)$. We have $t + t = 2t + 0$, and so $\psi(2t) = 2\psi(t) - \psi(0)$. By an easy induction we have $\psi(p_kt) = p_k\psi(t) - (p_k - 1)\psi(0)$. However, since $p_k t = 0$, this implies that $p_k(\psi(t) - \psi(0)) = 0$. But $\psi$ is an isomorphism, and so $\psi(t) \neq \psi(0)$, from which it follows that $G'$ does indeed have an element of order $p_k$. Similarly, $G'$ has elements of order $p_1,\dots,p_{k-1}$ and $Q$.\vs

\ni Now we have 
\[ n  =  |A|  =  Q(p_1 + \dots + p_k - k - 1) \; \geq \; QX^2/4\log X,\] whilst
\[ |A + A| \; \leq \; Q (p_1 + \dots + p_k)^2 \; \leq \; QX^4/(\log X)^2.\]
Thus our choice of $X$ guarantees that $|A + A| \leq Kn$. Now the claim proved above implies that if $G'$ contains a model for $A$ then $|G'| \geq Qp_1\dots p_k \geq Q2^X$. Since $n \leq QX^2/\log X$, we have
\[ \frac{|G'|}{n} \; \geq \; \frac{2^{X}\log X}{X^2} \; \geq \; e^{\frac{1}{6}\sqrt{K \log K}},\] as claimed.\endproof

\section{Acknowledgement} \ni The authors would like to thank Terry Tao for suggesting the term ``coset progression'', and for helpful remarks concerning the material of \S \ref{sec4}.


\begin{thebibliography}{99}

\bibitem{bilu} Y.~Bilu, \emph{Structure of sets with small sumset,} in Structure theory of set addition, Ast\'erisque \textbf{258} (1999), 77--108.

\bibitem{bogolyubov} N.~N.~Bogolyubov, \emph{Sur quelques propri\'et\'es arithm\'etiques des presque-p\'eriodes}, Ann. Chaire Math. Phys. Kiev \textbf{4} (1939), 185--194.

\bibitem{cassels}
J.~W.~S.~Cassels, \emph{An introduction to the geometry of numbers,} Corrected reprint of the 1971 edition. Classics in Mathematics. Springer-Verlag, Berlin, 1997. viii+344 pp (originally published in 1959).

\bibitem{chang}
M.~C.~Chang, \emph{A polynomial bound in Freiman's theorem},
Duke Math. J. \textbf{113} (2002), no. 3, 399--419.

\bibitem{dhp} J.~-M.~Deshouillers, F.~Hennecart and A.~Plagne, \emph{On small sumsets in $(\Z/2\Z)^n$,} Combinatorica \textbf{24} (2004), 53--68.

\bibitem{freiman} 
G.~Freiman, \emph{Foundations of a structural theory of set addition}, Translations of Mathematical Monographs \textbf{37}, Amer. Math. Soc., Providence, RI, USA, 1973. 

\bibitem{green-edinburgh} B.~J.~Green, \emph{Edinburgh-MIT lecture notes on Freiman's theorem,} in preparation for Online J. Analytic Combinatorics.

\bibitem{green-ruzsa} B.~J.~Green and I.~Z.~Ruzsa, \emph{Sets with small sumset and rectification,} to appear in Bull. London Math. Soc. 

\bibitem{plunnecke}
H.~Pl\"unnecke, \emph{Eigenschaften un Absch\"atzungen von Wirkingsfunktionen}, BMwF-GMD-22 Gesellschaft f\"ur Mathematik und Datenverarbeitung, Bonn (1969).

\bibitem{rudin} W.~Rudin, \emph{Fourier analysis on groups,} Interscience Tracts in Pure and Applied Mathematics, No. 12 Interscience Publishers (a division of John Wiley and Sons), New York-London 1962, 285 pp. 

\bibitem{ruzsa-sumdiff} I.~Z.~Ruzsa, \emph{On the cardinality of $A+A$ and $A-A$,} 
Combinatorics (Proc. Fifth Hungarian Colloq., Keszthely, 1976), Vol. II, pp. 933--938, 
Colloq. Math. Soc. J\'anos Bolyai \textbf{18} 
North-Holland, Amsterdam-New York, 1978. 

\bibitem{ruzsa-plun}
I.~Z.~Ruzsa, \emph{An application of graph theory to additive number theory}, Scientia, Ser. A. \textbf{3} (1989), 97--109.

\bibitem{ruzsa-freiman}
I.~Z.~Ruzsa, \emph{Generalized arithmetical progressions and sumsets}, 
Acta Math. Hungar. \textbf{65} (1994), no. 4, 379--388.

\bibitem{ruzsa-freiman-torsion} I.~Z.~Ruzsa, \emph{An analog of Freiman's theorem in groups}, Structure theory of set addition, Ast\'erisque \textbf{258} (1999), 323--326.


\end{thebibliography}
     \end{document}